\documentclass[reqno,11pt]{amsart}
\usepackage{graphicx} 
\usepackage{tikz-cd}
\usepackage{hyperref}
\usepackage{soul}
\usepackage{mathrsfs}  

\usepackage{enumitem}
\usepackage{marginnote}
\usepackage[nobysame]{amsrefs}
\usepackage{hyperref}
\usepackage[scr=rsfso]{mathalfa}
\usepackage{xcolor}
\usepackage{slashed}
\usepackage{subcaption}
\usepackage{mathtools}

\usepackage[nobysame]{amsrefs}
\usepackage{comment}
\usepackage{soul}
\usepackage{tikz}
\usepackage{marginnote}
\usepackage{anysize}
\usepackage{bbm}
\usepackage{cancel}
\usepackage[normalem]{ulem}
\usepackage[percent]{overpic} 

\theoremstyle{plain}
\newtheorem{theorem}[subsection]{Theorem}
\newtheorem{lemma}[subsection]{Lemma}

\newtheorem{proposition}[subsection]{Proposition}

\theoremstyle{definition}
\newtheorem{definition}[subsection]{Definition}

\theoremstyle{remark}

\newcommand{\RR}{\mathbb{R}}

\newcommand{\upper}{\uppercase\expandafter}

\newcommand{\supp}{\operatorname{supp}}

\newcommand{\calH}{\mathscr{H}}

\newcommand{\rank}{\operatorname{rank}}

\newcommand{\oM}{\overline{M}}\newcommand{\oV}{\overline{V}}\newcommand{\oH}{\overline{H}}

\begin{document}
\title{The semi-classical Weyl law on complete manifolds}
\author{Maxim Braverman} 
\address{Department of Mathematics,
        Northeastern University,
        Boston, MA 02115,
        USA
         }
\email{m.braverman@northeastern.edu}

\begin{abstract}
We prove that the semi-classical Schr\"odinger operator with growing potential on a complete Riemannian manifold satisfies the Weyl law.
\end{abstract}


\maketitle


H. Weyl, \cite{Weyl1911asymptotische},  established a remarkable asymptotic of the number of large eigenvalues of a Dirichlet Laplacian on a bounded domain in $\RR^n$, which is now called the Weyl law. This result was generalized and refined in many ways, see \cite{Ivrii2016weyl} for a review. In particular, a Weyl law for a Schr\"odinger operator on compact manifolds was established. Another version is the semi-classical Weyl law for the semi-classical Schr\"odinger operator on a compact manifold (see Section~\ref{SS:semiclassicalWeyl} for details). 

The Weyl law was established for Schr\"odinger operators on $\RR^n$ with some conditions on the growing rate of the potential, \cites{Rozenbljum74,tachizawa1992eigenvalue,Levendorskii96}. There are also some results for asymptotically hyperbolic manifolds \cites{Bonthonneau2015WeylLF, Moroianu2008}, asymptotically Euclidean manifolds \cite{CoriascoDoll2021}, and manifolds with cylindrical ends \cite{ManicciaPanarese2002}. 

In a recent paper, \cite{DaiYan2025heat}, X.~Dai and J.~Yan established both classical and semi-classical versions of the Weyl law for Schr\"odinger operators on manifolds of bounded geometry with some conditions on the growth of the potential. Their method was to obtain an asymptotic expansion of the heat kernel, which is a much stronger result than the Weyl law. 

In this note, we establish the semi-classical Weyl law on an arbitrary complete Riemannian manifold without any restrictions on the potential (except that it grows at infinity). We achieve it by comparing the semi-classical spectral counting functions on the complete manifold and on a certain compact manifold. This is done using the IMS localization formula \cite{CFKS}*{\S3.1} by a method similar to what we used in \cite{BrFar2}*{\S4} and \cite{Br-cob}.

We should note that the classical Weyl law in non-compact setting, \cite{DaiYan2025heat}*{Th.~1.8}, is more subtle and is not approachable by our methods.  

\subsection*{Acknowledgment} I am grateful to the Max Planck Institute for Mathematics for its hospitality and for providing an excellent research environment. I also thank Junrong Yan for valuable discussions.
\subsection{The semi-classical Weyl law}\label{SS:semiclassicalWeyl}
Let $M$ be a complete $n$-dimensional manifold.  We denote by $\Delta$ the positive Laplacian on $M$. Let $V:M\to [0,\infty)$ be a smooth function which tends to plus infinity as $x\to\infty$. Let 
\begin{equation}\label{E:Schrodinger}
    H_{\hbar} \ := \ \hbar^2\Delta\ + \ V,
\end{equation}
be the {\em semi-classical Schr\"odinger operators on $M$}. Then $H_\hbar$ is essentially self-adjoint, non-negative, and has a discrete spectrum.

For a self-adjoint operator $A$ with discrete spectrum,  we define its {\em spectral counting function} $N(A,\lambda)$ as the number of eigenvalues of $A$ which are smaller than $\lambda$. In general,  $N(A,\lambda)$ can be infinite, but our assumptions on $H_\hbar$ guarantee that its spectral counting function is finite. 

Our main result is the following

\begin{theorem}\label{T:semiclassicaWeyl}
For any $\lambda\in \RR$, we have 
\begin{equation}\label{E:semiclassicaWeyl}
    {(2\pi\hbar)}^n\, N(H_\hbar,\lambda)\ = \  
    {\rm Vol}\big\{\, (x,\xi)\in T^*M:\, |\xi|^2+V(x)\le \lambda\,\big\}\ + \ o(1),
\end{equation}
as $\hbar\to 0$.
\end{theorem}

The rest of the paper is devoted to the proof of this theorem, which is based on comparing the eigenvalue counting function on $M$ and on a certain compact manifold.

\subsection{The relative semi-classical Weyl law}\label{SS:ralativeWeyl}
Let $M_1$ and $M_2$ be complete (possibly compact) $n$-dimensional manifolds.
Let $V_j\in [0,\infty)$ be smooth functions. If $M_j$ is non-compact, we assume that $V(x)$ tends to plus infinity as $x\to\infty$. The equation \eqref{E:Schrodinger} defines the semi-classical Schr\"odinger operators  $H_{\hbar,j}$ on $M_j$.

\begin{definition}
Fix a real number $\lambda>0$. We say that the pairs $(M_1,V_1)$ and $(M_2,V_2)$ are {\em $\lambda$-equivalent} if there exist a relatively compact open subsets $U_j\subset M_j$ and an isometric diffeomorphism $\Phi:U_1\to U_2$ such that 
\begin{enumerate}
    \item $ V_j^{-1}\big([0,\lambda]\big)\subset U_j$;
    \item \(V_2\big(\Phi(x)\big) \ = \  V_1(x)\), for all $x\in U_1$.
\end{enumerate}
\end{definition}

We have the following {\em relative} version of the semi-classical Weyl's law:

\begin{proposition}\label{P:relativeWeylsc}
If $(M_1,V_1)$ and $(M_2,V_2)$ are $\lambda$–equivalent, then there exists a constant $c>0$ (which depends on $\lambda$ and all the other data, but is independent of $\hbar$) such that 
\begin{equation}\label{E:relativeWeylsc}
    N(H_{\hbar,1},\lambda)  \ \ge \ 
    N(H_{\hbar,2},\lambda-c\hbar^2).
\end{equation}
\end{proposition}

Before proving the proposition, we show how it implies the semi-classical Weyl law.

\subsection{The proof of the semi-classical Weyl law}\label{SS:pr}
Fix $\lambda$ and and an open relatively compact set $U\subset M$ which contains the set $\big\{\, (x,\xi)\in T^*M:\, |\xi|^2+V(x)\le \lambda\,\big\}$. Let $\oM$ be any compact Riemannian manifold that contains $U$. Let $\oV:\oM\to \RR$ be a smooth function whose restriction to $U$ equals $V$ and whose restriction to $\oM\backslash U$ is bigger than $\lambda$. Then $(M,V)$ and $(\oM,\oV)$ are $\lambda$-equivalent.  We denote by $\oH_\hbar$ the semi-classical Schr\"odinger operator on $\oM$:
\[
    \oH_\hbar\ := \ \hbar^2\,\Delta + \oV.
\]

It is a classical result that on compact manifolds, the eigenvalue counting function satisfies the semi-classical Weyl law, cf. \cite{Zworski12book}*{Th.~6.8}. It follows that 
\begin{multline}\label{E:Ncompact-}
        (2\pi\hbar)^n\,N(\oH_\hbar,\lambda)\ = \ {\rm Vol}\big\{\, (x,\xi)\in T^*M:\, |\xi|^2+V(x)\le \lambda\,\big\} +o(1)
    \\ = \ (2\pi\hbar)^n\,N_h(\oV,\lambda-c\hbar^2)+o(1).
\end{multline}

Setting $M_1=M$, $M_2= \oM$ in \eqref{E:relativeWeylsc} and using \eqref{E:Ncompact-} we now obtain
\begin{equation}\label{E:M>oM}
    (2\pi\hbar)^n\,N(H_\hbar,\lambda)\ \ge \ (2\pi\hbar)^n\,N(\oH_\hbar,\lambda-c\hbar^2)
    \ = \ {\rm Vol}\big\{\, (x,\xi)\in T^*M:\, |\xi|^2+V(x)\le \lambda\,\big\} +o(1).
\end{equation}

For any $\epsilon >0$ there exists $\delta>0$ such that $V^{-1}\big([0,\lambda+\delta]\big)\subset U$ and
\begin{equation}\label{E:delta-epsilon}
    {\rm Vol}\big\{\, (x,\xi)\in T^*M:\, \lambda\le |\xi|^2+V(x)\le \lambda+\delta\,\big\}
    \ <\ \epsilon. 
\end{equation}
Setting  $M_1=\oM$ and $M_2=M$ and replacing $\lambda$ with $\lambda+\delta$ in Proposition~\ref{P:relativeWeylsc}, we conclude that there exists a constant $c'>0$ such that 
\[
    N(\oH_\hbar,\lambda+\delta) \ \ge \ N(H_\hbar,\lambda+\delta-c'\hbar^2).
\]
Hence, using \eqref{E:Ncompact-} and \eqref{E:delta-epsilon}, we have
\begin{equation}\label{E:oM>M}
    N(\oH_\hbar,\lambda)+ \frac{\epsilon}{(2\pi \hbar)^n} \ > \ N(\oH_\hbar,\lambda+\delta) \ \ge \ 
    N(H_\hbar,\lambda+\delta-c'\hbar^2).
\end{equation}
If $c'\hbar^2< \delta$ then $N(H_\hbar,\lambda+\delta-c'\hbar^2)\ge N(H_\hbar,\lambda)$. Thus \eqref{E:oM>M} implies that for every $\epsilon>0$ there exists $\hbar_\epsilon>0$ such that for all $\hbar<\epsilon$ we have 
\begin{equation}\notag
    (2\pi\hbar)^n\,N(H_\hbar,\lambda)\ \le \ (2\pi\hbar)^n\,N(\oH_\hbar,\lambda) + \epsilon.
\end{equation}
It follows that 
\begin{equation}\notag
    (2\pi\hbar)^n\,N(H_\hbar,\lambda)\ \le \ (2\pi\hbar)^n\,N(\oH_\hbar,\lambda) + o(1)
    \ = \ {\rm Vol}\big\{\, (x,\xi)\in T^*M:\, |\xi|^2+V(x)\le \lambda\,\big\} +o(1).
\end{equation}
Together with \eqref{E:M>oM} this proves the theorem. 
\hfill$\square$

To finish the proof of Theorem~\ref{T:semiclassicaWeyl} it remains to prove Proposition~\ref{P:relativeWeylsc}. We do it using a very useful formula called Ismagilov–Morgan–Simon (IMS) localization,   cf. \cite{CFKS}*{\S3.1}.

\subsection{Ismagilov–Morgan–Simon (IMS) localization formula}\label{SS:IMS}
We use the following IMS localization formula, cf. \cite{CFKS}*{\S3.1}).
\begin{lemma}\label{L:IMS} 
Suppose $\phi$ and $\psi$ are smooth functions with $\psi^2+\phi^2=1$.
Then the following operator identity holds
  \begin{equation}\label{E:IMS}
        H=\phi H\phi+\psi H \psi+\frac12[\phi,[\phi,H]]+\frac12[\psi,[\psi,H]].
  \end{equation}
\end{lemma}
The proof is given in \cite{CFKS}*{\S3.1} (see also \cite{BrFar2}*{Lemma~4.10} or \cite{BrSil}*{Lemma~8.2}).
Note that if $H$ is the semi-classical Schr\"odinger operator \eqref{E:Schrodinger}, then 
\begin{equation}\label{E:phphiH}
    \frac12\,[\phi,[\phi,H]]  \ = \ -\hbar^2\,|d\phi|^2, \qquad 
    \frac12\,[\psi,[\psi,H]]  \ = \ -\hbar^2\,|d\psi|^2.
\end{equation} 

\subsection{Estimates on $N(H_{\hbar,j},\lambda)$}\label{SS:estimates}
We use the notation of Section~\ref{SS:ralativeWeyl}.
For $j=1,2$, choose functions $\phi_j,\,\psi_j:M_j\to[0,1]$ such that
\begin{enumerate}
\item $\phi_j^2+\psi_j^2=1$;
\item $\supp \phi_j\in U_j$;
\item $\supp \psi_j \cap V_j^{-1}\big([0,\lambda]\big)= \emptyset$,
\item $\phi_1(x)= \phi_2(\Phi(x))$ for every $x\in U_1$.
\end{enumerate}
Then
\begin{equation}\label{E:psiHpsi>}
    \psi_j\,H_{\hbar,j}\,\psi_j\  > \ \psi_j^2\lambda.
\end{equation}
From \eqref{E:phphiH} we conclude that there exists a constant $c>0$ such that 
\begin{equation}\label{E:phphiH<}
    \big|\frac12\,[\phi_j,[\phi_j,H_{\hbar,j}]]\big| \ < \ \frac{c}2\hbar^2,\qquad
    \big|\frac12\,[\psi_j,[\psi_j,H_{\hbar,j}]]\big|  \ < \ \frac{c}2\hbar^2.
\end{equation}

Finally, let $P_{\hbar,j}$ denote the orthogonal projection onto the span of the eigenfunctions of $H_{\hbar,j}$ with eigenvalues less than $\lambda$. By definition, $\rank(P_{\hbar,j})= N(H_{\hbar,j},\lambda)$. Then 
\begin{equation}\label{E:H+P>lambda}
    H_{\hbar,j}+ \lambda P_{\hbar,j}\ \ge \ \lambda \qquad \Longrightarrow\qquad
    \phi_1\,\big(\,H_{\hbar,j}+ 
    \lambda P_{\hbar,j}\,\big)\,\phi_1\ \ge \ \lambda\,\phi_1^2.
\end{equation}

Using the diffeomorphism $\Phi:U_1\to U_2$ we can view $\phi_1 P_{\hbar,1}\phi_1$ as an operator acting on functions on $M_2$. With a slight abuse of notation, we occasionally denote the same operator by  $\phi_2P_{\hbar,1}\phi_2$.
Then the operator 
\[
    \phi_2 H_{\hbar,2}\phi_2 + \lambda\phi_2 P_{\hbar,1}\phi_2
\]
can be identified with $ \phi_1\big(\,H_{\hbar,1}+ \lambda P_{\hbar,1}\big)\,\phi_1$. It now follows from 
\eqref{E:H+P>lambda} that 
\begin{equation}\label{E:phiH+B>lambda}
    \phi_2 H_{\hbar,2}\phi_2 + 
    \lambda\phi_2 P_{\hbar,1}\phi_2\ \ \ge \ \lambda\phi_2^2.
\end{equation}
Combining \eqref{E:IMS}, \eqref{E:psiHpsi>}, \eqref{E:phphiH}, and \eqref{E:phiH+B>lambda}, we obtain
\begin{equation}\label{E:H2>}
    H_{\hbar,2}+\phi_2 P_{\hbar,1}\phi_2\ \ge \ \lambda-c\hbar^2.
\end{equation}

Proposition~\ref{P:relativeWeylsc} follows now from the following general lemma (cf. \cite{ReSi4}*{p. 270}):
\begin{lemma}\label{L:general}
   Assume that $A, B$ are self-adjoint operators in a Hilbert space
   $\calH$ such that \linebreak $\rank(B)\le k$ and there exists $\mu>0$
   such that
   \[
        \langle (A+B)u,u\rangle \ge \mu\langle u,u\rangle
                \quad \text{for any} \quad u\in\operatorname{Dom}(A).
   \]
   Then $N(A,\mu)\le k$.
\end{lemma}
To finish the proof of the proposition, we simply set $B:= \lambda \phi_2 P_{\hbar,1}\phi_2$ and observe  that 
\[
    \rank\big(\phi_2 P_{\hbar,1}\phi_2)\ \le \ \rank(P_{\hbar,1})\  = \ N(H_{\hbar,1},\lambda)
\]
and set $\mu= \lambda-c\hbar^2$ in Lemma~\ref{L:general}. \hfill $\square$



\begin{bibdiv}
\begin{biblist}

\bib{Bonthonneau2015WeylLF}{article}{
      author={Bonthonneau, Yannick~Guedes},
       title={Weyl laws for manifolds with hyperbolic cusps},
        date={2015},
     journal={arXiv: Spectral Theory},
         url={https://api.semanticscholar.org/CorpusID:119160723},
}

\bib{Br-cob}{article}{
      author={Braverman, Maxim},
       title={New proof of the cobordism invariance of the index},
        date={2002},
     journal={Proc. Amer. Math. Soc.},
      volume={130},
      number={4},
       pages={1095\ndash 1101},
}

\bib{BrFar2}{article}{
      author={Braverman, Maxim},
      author={Farber, Michael},
       title={Novikov type inequalities for differential forms with
  non-isolated zeros},
        date={1997},
     journal={Math. Proc. Cambridge Philos. Soc.},
      volume={122},
       pages={357\ndash 375},
}

\bib{BrSil}{article}{
      author={Braverman, Maxim},
      author={Silantyev, Valentin},
       title={Kirwan-{N}ovikov inequalities on a manifold with boundary},
        date={2006},
     journal={Trans. Amer. Math. Soc.},
      volume={358},
       pages={3329\ndash 3361},
}

\bib{CoriascoDoll2021}{article}{
      author={Coriasco, Sandro},
      author={Doll, Moritz},
       title={Weyl law on asymptotically {E}uclidean manifolds},
        date={2021},
     journal={Annales Henri Poincar{\'e}},
      volume={22},
      number={2},
       pages={447\ndash 486},
}

\bib{CFKS}{book}{
      author={Cycon, H.L.},
      author={Froese, R.G.},
      author={Kirsch, W.},
      author={Simon, B.},
       title={Schr\"odinger operators with applications to quantum mechanics
  and global geometry},
      series={Texts and Monographs in Physics},
   publisher={Springer-Verlag},
        date={1987},
}

\bib{DaiYan2025heat}{article}{
      author={Dai, Xianzhe},
      author={Yan, Junrong},
       title={Heat kernel approach to the {W}eyl law for {S}chr{\"o}dinger
  operators on non-compact manifolds},
        date={2025},
     journal={arXiv preprint arXiv:2504.15551},
         url={https://arxiv.org/abs/2504.15551},
        note={Submitted on 22 Apr 2025},
}

\bib{Ivrii2016weyl}{article}{
      author={Ivrii, Victor},
       title={100 years of {W}eyl's law},
        date={2016},
     journal={Bulletin of Mathematical Sciences},
      volume={6},
      number={3},
       pages={379\ndash 452},
         url={https://doi.org/10.1007/s13373-016-0089-y},
}

\bib{Levendorskii96}{article}{
      author={Levendorski\u{\i}, S.~Z.},
       title={Spectral asymptotics with a remainder estimate for
  {S}chr\"{o}dinger operators with slowly growing potentials},
        date={1996},
        ISSN={0308-2105},
     journal={Proc. Roy. Soc. Edinburgh Sect. A},
      volume={126},
      number={4},
       pages={829\ndash 836},
         url={https://doi-org.ezproxy.neu.edu/10.1017/S030821050002309X},
      review={\MR{1405759}},
}

\bib{ManicciaPanarese2002}{article}{
      author={Maniccia, Luciano},
      author={Panarese, Paola},
       title={Eigenvalue asymptotics for a class of md-elliptic
  pseudo-differential operators on manifolds with cylindrical exits},
        date={2002},
     journal={Annali di Matematica Pura ed Applicata},
      volume={181},
      number={3},
       pages={283\ndash 308},
}

\bib{Moroianu2008}{article}{
      author={Moroianu, Sergiu},
       title={Weyl laws on open manifolds},
        date={2008},
     journal={Mathematische Annalen},
      volume={340},
       pages={1\ndash 21},
}

\bib{ReSi4}{book}{
      author={Reed, M.},
      author={Simon, B.},
       title={Methods of modern mathematical physics {I}{V}: {Analysis} of
  operators},
   publisher={Academic Press},
     address={London},
        date={1978},
}

\bib{Rozenbljum74}{article}{
      author={Rozenbljum, G.~V.},
       title={Asymptotic behavior of the eigenvalues of the {S}chr\"{o}dinger
  operator},
        date={1974},
     journal={Mat. Sb. (N.S.)},
      volume={93 (135)},
       pages={347\ndash 367, 487},
      review={\MR{0361470}},
}

\bib{tachizawa1992eigenvalue}{article}{
      author={Tachizawa, Kazuya},
       title={Eigenvalue asymptotics of {S}chr{\"o}dinger operators with only
  discrete spectrum},
        date={1992},
     journal={Publications of the Research Institute for Mathematical
  Sciences},
      volume={28},
      number={6},
       pages={943\ndash 981},
}

\bib{Weyl1911asymptotische}{article}{
      author={Weyl, Herman},
       title={{\"U}ber die asymptotische verteilung der eigenwerte},
        date={1911},
     journal={Nachrichten von der Gesellschaft der Wissenschaften zu
  G{\"o}ttingen, Mathematisch-Physikalische Klasse},
       pages={110\ndash 117},
}

\bib{Zworski12book}{book}{
      author={Zworski, Maciej},
       title={Semiclassical analysis},
      series={Graduate Studies in Mathematics},
   publisher={American Mathematical Society, Providence, RI},
        date={2012},
      volume={138},
        ISBN={978-0-8218-8320-4},
         url={https://doi-org.ezproxy.neu.edu/10.1090/gsm/138},
      review={\MR{2952218}},
}

\end{biblist}
\end{bibdiv}
\end{document}